\documentclass[11pt]{article}

\usepackage[mathscr]{eucal}
\usepackage{color}
\usepackage{amsmath}
\usepackage{amssymb}
\usepackage{amsfonts}
\usepackage{amscd}
\usepackage{amsthm}
\usepackage{latexsym}
\usepackage{graphicx}
\usepackage{subfig}

\def\be{\begin{equation}}
\def\ee{\end{equation}}
\def\bse{\begin{subequations}}
\def\ese{\end{subequations}}

\usepackage{latexsym}

\let\er\eqref

\let\be\beta

\newcommand{\R}{{\mathbb R}}

\newtheorem{theorem}{Theorem}
\newtheorem{lemma}[theorem]{Lemma}
\newtheorem{proposition}[theorem]{Proposition}

\def\bse{\begin{subequations}}
\def\ese{\end{subequations}}

\usepackage{geometry}

\title{Global existence in the critical and subcritical cases to the Fisher-KPP model with nonlocal nonlinear reaction}

\author{ Shen Bian\footnote{Beijing University of Chemical Technology, 100029, Beijing. Email: \texttt{bianshen66@163.com}. Partially supported by National Science Foundation of China (Grant No. 11501025).} }
\date{}

\begin{document}
\let\cleardoublepage\clearpage%svinei kenes selides

\maketitle

\begin{abstract}
The Cauchy problem considered in this paper is the following
\begin{align}\label{nkpp0}
  \left\{
    \begin{array}{ll}
      u_t=\Delta u+u^\alpha\left(M_0- \int_{\R^n} u(x,t)dx\right),\quad   & x \in \R^n, t>0,
  \\
   u(x,0)=U_0(x)\geq 0,\quad & x \in \R^n.
    \end{array}
  \right.
  \end{align}
  where $M_0>0, \alpha >1, n \ge 3$. When the coefficient $M_0-\int_{\R^n} u(x,t) dx$ remains positive, \er{nkpp0} is analogous to
  \begin{align}\label{fujita}
  \left\{
    \begin{array}{ll}
      u_t=\Delta u+u^\alpha,\quad   & x \in \R^n, t>0,
  \\
   u(x,0)=U_0(x)\geq 0,\quad & x \in \R^n.
    \end{array}
  \right.
  \end{align}
It is well known that when $1<\alpha \le 1+2/n$, the local solution of \er{fujita} blows up in finite time as long as the initial value is nontrivial. The present paper forms a contrast to \er{fujita} and shows the global existence of solutions to \er{nkpp0} for $1<\alpha \le 1+2/n$ by dealing with the mathematical challenge which is from the nonlocal term $\int_{\R^n} u dx$. It's proved that when $1<\alpha <1+2/n$, such a global bound is obtained for any positive $M_0$ and any non-negative initial data. While if $\alpha=1+2/n$, then the global solution does exist for sufficiently small $M_0$ and any non-negative initial data. Furthermore, the large time behavior of the global solution is also discussed for $\alpha=1+2/n$. Besides, this paper establishes the hyper-contractivity of a global solution in $L^\infty(\R^n)$ with $U_0 \in L^1(\R^n)$ for the case $\alpha=1+2/n$.
\end{abstract}
%%%%%%%%%%%%%%%%%%%%%%%%%%%%%%%%%%%%%%%%%%%%%%%%%%%%%%%%%%%%%%%%%%%%%%%%%%%%%%%%%%%%%%%%%%%%%%%%%%%%%%%%%%%%%%%%%%%%%%%%%%%%%%%%%%%%%%%%%%%%%%%%%%%%%%%%%%%%%%%%%%%%%%%%%%%%%%%%%%%%%%%%%%%%%%%%%%%%%%%%%%%%%%%%%%%%%%%%%%%%%%%%%%%%%%%%%%%%%%%%%%%%%%%%%%%%%%%%%%%%%%%%%%%%%%%%%%%%%%%%%%%%%%%%%%%%%%%%%%%%%%%%%%%%%%%%%%%%%%%%%%%%%%%%%%%%%%%%%%%%%%%%%%%%%%%%%%%%%%%%%%%%%%%%%%%%%%%%%%%%%%%%%%%%%%%%%%%%%%%%%%%%%%%%%%%%%%%%%%%%%%%%%%%%%%%%%%%%%%%%%%%%%%%%%%%%%%%%%%%%%%%%%%%%%%%%%%%%%%%%%%%%%%%%%%%%%%%%%%%%%%%%%%%%%%%%%%%%%%%%%%%%%%%%%%%%%%%%%%%%%%%%%%%%%%%%%%%%%%%%%%%%%%%%%%%%%%%%%%%%%%%%%%%%%%%%%%%%%%%%%%%%%%%%%%%%%%%%%%%%%%%%%%%%%%%%%%%%%%%%%%%%%%%%%%%%%%%%%%%%%%%%%%%%%%%%%%%%%%%%%%%%%%%%%%%%%%%%%%%%%%%%%%%%%%%%%%%%%%%%%%%%%%%%%%%%%%%%%%%%%%%%%%%%%%%%%

\section{Introduction}\label{sec1}
\def\theequation{1.\arabic{equation}}\makeatother
\setcounter{equation}{0}
\def\thetheorem{1.\arabic{theorem}}\makeatother
\setcounter{theorem}{0}

This paper deals with the following semilinear equation with nonlocal reaction term
\begin{align}\label{nkpp}
\left\{
  \begin{array}{ll}
    u_t=\Delta u+u^\alpha\left(M_0- \int_{\R^n} u(x,t)dx\right),\quad   & x \in \R^n, t>0,
\\
 u(x,0)=U_0(x) \geq 0,\quad &x \in \R^n,
  \end{array}
\right.
\end{align}
where $M_0>0, \alpha>1, n \ge 3$. \er{nkpp} is related to many equations arising from population dynamics and combustion theory with Fisher-KPP type reaction terms \cite{f37,kpp}. Here $u$ is the density of the population, $M_0$ can be viewed as a measure of the strength of the reaction mechanism and $M_0 u^\alpha$ induces a power-like growth for low density populations. The nonlocal term $-\int_{\R^n} u(x,t)dx$ describes the influence of the total mass on the growth of the population \cite{NT13} which counteracts the blow-up tendency produced by growth factor. Since they present densities, solutions to \er{nkpp} which are biologically meaningful must satisfy
\begin{align}
u(x,t) \ge 0,\quad t>0.
\end{align}
This can also be verified by the maximum principle since the non-negativity of the initial data. In view of this, it's reasonable to require throughout that the initial data $U_0\in C(\R^n)$ decaying at infinity be non-negative. We shall deal with the solution $u(x,t)$ which satisfies the integral equation in spirit of \cite{fuj1,ln92}
\begin{align}\label{property}
  \begin{array}{ll}
    u(x,t)=&\frac{1}{(4\pi t)^{n/2}} \int_{\R^n} e^{-\frac{|x-y|^2}{4t}} U_0(y) dy \\
& +\int_0^t \int_{\R^n} \frac{1}{(4 \pi (t-s))^{n/2}}e^{-\frac{|x-y|^2}{4t(t-s)}} u^\alpha (y,s) \left(M_0-\int_{\R^n} u(z,s)dz \right) dy ds.
  \end{array}
\end{align}
We will give a detailed introduction to \er{property} in Section \ref{subsec2}.

Before we turn to the study of \er{nkpp} in the direction concerned in this paper, we mention that diffusion equations with nonlocal reactions in bounded domain
\begin{align}\label{expononlocal}
u_t=\Delta u+F(u)
\end{align}
have been considered by a number of authors \cite{bb82,BCL15,hy95,lcl09,Pao92,VVpp,ww11,ww96,Yi96}. Bebernes and Bressan \cite{bb82} studied an ignition model for a compressible reactive gas and Pao \cite{Pao92} discussed combustion theory which is a reaction-diffusion equation with nonlocal exponential reaction term
\begin{align}\label{expononlocal}
\left\{
  \begin{array}{ll}
    u_t=\Delta u+\sigma \left(e^{\gamma u} +b \int_{\Omega}e^{\gamma u} dx \right),& x \in \Omega, t>0, \\
    u(x,t)=0,& x \in \partial \Omega, t>0, \\
    u(x,0)=U_0(x),& x \in \Omega.
  \end{array}
\right.
\end{align}
There exists a bounded value $\sigma^\ast$ related to the boundedness of $|\Omega|$ and the first eigenvalue of $-\Delta$ with homogeneous Dirichlet boundary condition such that for $\sigma<\sigma^\ast$, there is a unique global solution to \er{expononlocal}.

As the parabolic inverse problem, Hu and Yin \cite{hy95} studied the following problem with nonlocal dampening term
\begin{align}\label{massnonlocal}
\left\{
  \begin{array}{ll}
    u_t=\Delta u+u^p-\frac{1}{|\Omega|} \int_{\Omega} u^p dx, & x \in \Omega, t>0, p>1, \\
\frac{\partial u}{\partial \nu}=0, & x \in \partial \Omega, t>0, \\
u(x,0)=U_0(x),& x \in \Omega.
  \end{array}
\right.
\end{align}
This typical structure admits mass conservation and the corresponding Lyapunov functional $J(t)=\frac{1}{2}\int_{\Omega} |\nabla u|^2dx+\frac{1}{p+1}\int_{\Omega} u^{p+1} dx$ is non-increasing in time. In view of the boundedness of $|\Omega|$, they proved that the solution of \er{massnonlocal} blows up if $J(0)$ is suitably large by convexity arguments. Later, Wang and Wang \cite{ww96} considered the local dampening reaction term in the following problem
\begin{align}\label{localdamp}
\left\{
  \begin{array}{ll}
    u_t=\Delta u+\int_{\Omega} u^q dx-k u^p(t,x), & x \in \Omega, ~t>0,~ p,q \ge 1,~k>0, \\
u(x,0)=U_0(x),& x \in \Omega.
  \end{array}
\right.
\end{align}
Subject to Neumann boundary condition $\frac{\partial u}{\partial \nu}=0, x \in \partial \Omega, t>0.$ By comparing with solutions of an ODE they proved that for $q<p,$ \er{localdamp} has a global solution; if $p=q=1$ or $p=q>1$ and $|\Omega|\le k$, then the solution of problem \er{localdamp} blows up in finite time provided that $U_0(x)$ is not identically zero; if $q>p$, then \er{localdamp} has a global solution when $U_0(x)\le f^\ast$ and the solution blows up in finite time when $U_0(x)>f^\ast$ where $f^\ast$ is a constant depending on $k$ and $|\Omega|.$

Nonlocal type reaction terms can also describe Darwinian evolution of a structured population density or the behavior of cancer cells with therapy as well as chemotherapy, we refer the interested reader to the models found in \cite{Lorz:2011hl,Lorz:2013vp,vol1,vol2}.

Recently, Tello and Negreanu \cite{NT13} considered the logistic growth factor involving the coexistence of local and nonlocal consumption of resources
\begin{align}\label{nt13}
u_t=\Delta u+\chi u^2+u\left(a_0-a_1 u-\frac{a_2}{|\Omega|}\int_{\Omega} u dx\right),~~\chi,a_0,a_1,a_2>0
\end{align}
with Neumann boundary condition $\frac{\partial u}{\partial \nu}=0, x\in \partial \Omega$. Following a comparison argument they showed that if $a_1>2\chi+a_2$, then the solution of \er{nt13} converges to $\frac{a_0}{a_1+a_2}$ as time tends to infinity. Actually, as population grows, the competitive effect of the local term $a_1 u$ becomes more influential than the nonlocal term $a_2\int_{\Omega} u dx$ and the effect of the total mass can be ignored compared with the local term. As it was stated in \cite{NT13} that ``it seems to conjecture that the dampening effect of the nonlocal terms might lead to an even more effective homogenization". Therefore, the main aim of this paper is to explore the influence of the total mass on solutions of \er{nkpp}.

As for the remarkable difference between \er{nkpp} and \er{massnonlocal},\er{localdamp},\er{nt13}, we easily find that the estimates for well-posedness of solutions heavily rely on the boundedness of $|\Omega|$, which holds for \er{massnonlocal},\er{localdamp},\er{nt13} but not for \er{nkpp}. We remark that this comparison is more obvious between
\begin{align}\label{upbdd}
\left\{
  \begin{array}{ll}
     u_t=\Delta u+u^\alpha,\quad   & x \in \Omega, t>0, \alpha>1,
  \\
   u(x,0)=U_0(x) \geq 0,\quad & x \in \Omega,
  \end{array}
\right.
\end{align}
and
\begin{align}\label{fujita1}
\left\{
  \begin{array}{ll}
    u_t=\Delta u+u^\alpha,\quad   & x \in \R^n, t>0, \alpha>1,
  \\
  u(x,0)=U_0(x) \geq 0,\quad & x \in \R^n.
  \end{array}
\right.
\end{align}
We firstly mention that \er{upbdd} is working in a bounded domain $\Omega$. For \er{upbdd} with Dirichlet boundary condition $u\big|_{\partial \Omega}=0$ \cite{Per2015}, the first eigenvalue $\lambda_1$ of the operator $-\Delta$ is positive and associated with a positive eigenfunction $w_1(x)$. For $\alpha>1$, \er{upbdd} admits a global solution under the smallness condition $U_0(x)\le \left( \displaystyle \min_{\Omega} \frac{\lambda_1^{\frac{1}{\alpha-1}}}{w_1(\cdot)}\right) w_1(x)$. While the solution blows up in finite time for large initial data $\int_{\Omega} U_0(x) w_1(x) dx>\lambda_1^{1/(\alpha-1)} \int_{\Omega} w_1(x)dx$. On the other hand, nonexistence result was obtained for \er{upbdd} with Neumann boundary condition and no size condition on the initial data is required \cite{Per2015}. Now we turn to the whole space. The study of \er{fujita1} goes back to the fundamental work of Fujita \cite{fuj1}, he proved that if $1<\alpha<1+2/n,$ then \er{fujita1} has no non-negative global solution for any non-trivial initial data. The same is true when $\alpha=1+2/n$ as was proved by Hayakawa \cite{kh73} for $n=1,2$ and by Kobayashi, Sirao and Tanaka \cite{tana77} for general $n$. Especially for $\alpha=1+2/n$ with general $n$, Weissler \cite{weis81} demonstrated the nonexistence of a non-negative global solution by showing the unboundedness of solutions in $L^1(\R^n)$ which doesn't hold true in our model.

Now we are concerned \er{nkpp} that the
nonlocal reaction term $u^\alpha \int_{\R^n} u dx$ results in the lack of a good Lyapunov functional and the lack of comparison principle \cite{fuj1,weis81}. The main objective of the present one is to illustrate that the mortality or dampening term $-u^\alpha \int_{\R^n} u dx$ has a crucial influence on the solutions of the semilinear equation \er{nkpp} to exist for all time. To analyze this issue, we define the total mass
\begin{align}
m(t)=\int_{\R^n} u(x,t) dx
\end{align}
which satisfies
\begin{align}\label{mt}
\frac{d}{dt}m(t)=\left(M_0-m(t)\right) \int_{\R^n} u^\alpha dx
\end{align}
by integrating \er{nkpp} over $\R^n$. If the initial mass
\begin{align*}
m_0:=\int_{\R^n} U_0(x) dx>M_0,
\end{align*}
then we can see that $m(t)$ decreases in time and
\begin{align}
M_0 \le m(t) \le m_0
\end{align}
for all $t \ge 0.$ Thus we find
\begin{align*}
u_t=\Delta u+u^\alpha \left(M_0-m(t)\right) \le \Delta u
\end{align*}
which implies that $u(x,t)$ is a subsolution of the heat equation $v_t=\Delta v$ with the same initial data $U_0(x)$. Thus the solution of \er{nkpp} exists globally by the comparison principle. On the contrary, when the initial mass $m_0<M_0$, then $m(t)$ increases in time and
\begin{align} \label{m00}
m_0 \le m(t) \le M_0.
\end{align}
Therefore, without loss of generality, we assume that the initial mass satisfies
\begin{align}\label{m0xiaoyuM0}
m_0:=\int_{\R^n} U_0 dx<M_0
\end{align}
throughout this paper. In this sense, $M_0$ can be considered as the carrying capacity \cite{bp07}.

As pointed out in \er{fujita1} with the absence of nonlocal dampening term, all the solutions blow up in finite time for $1<\alpha\le 1+2/n.$ In the case of $m_0<M_0,$ we know that $M_0-m(t)$ is non-negative for all times which inspires us to have a natural guess that the problem \er{nkpp} might have no global solutions for $1<\alpha\le 1+2/n$. However, the present paper will give a negative answer to this observation.

Under the assumptions \er{m0xiaoyuM0}, an attempt to understand the global existence of solution $u(x,t)$ to \er{nkpp} is made for $1<\alpha \le 1+2/n$. Our main results can be summarized as follows:
\begin{itemize}
  \item $\alpha=1+2/n:$ Theorem \ref{critical} gives the existence and the decay property of a time global classical solution to \er{nkpp} when $U_0 \in C(\R^n) \cap L^\infty(\R^n)$ and $M_0$ is smaller than a uniform constant (see Lemma \ref{cstarexist}). In Theorem \ref{longcritical}, a precise large time behavior of the global classical solution $u(x,t)$ with the aid of \er{property} is obtained when global existence prevails. Indeed, it is not clear in the case $\alpha=1+2/n$ that for a given $U_0 \in L^1(\R^n)$, there exists a bounded global solution to \er{nkpp}. However, if $U_0 \in L^1 \cap L^\infty(\R^n)$, there is certainly a uniformly bounded global solution by Theorem \ref{critical}. Here we emphasize that the initial regularity $U_0 \in L^1 \cap L^\infty(\R^n)$ can be relaxed in deriving necessary conditions for the existence of a global solution to \er{nkpp}. Actually, to an initial data which is originally in $L^1(\R^n)$ but not in $L^\infty(\R^n)$, we associate a solution which at almost any time $t>0$ is in $L^k(\R^n)$ for $k$ arbitrarily large (hyper-contractivity) and furthermore in $L^\infty(\R^n)$ (ultra-contractivity), see Theorem \ref{hypercontractivity} for details.
  \item $1<\alpha<1+2/n:$ we prove that there exists a global classical solution without any restriction on $M_0$ where the initial data $U_0$ is a non-negative bounded continuous function, see Theorem \ref{subcritical}.
\end{itemize}

Our main results are stated and proved in Section \ref{sec3} and \ref{sec4}. We shall prepare several lemmas in Section \ref{sec2} which play an important role in Section \ref{sec3} and \ref{sec4}. Section \ref{sec5} contains some concluding remarks and guesses.

\section{Preliminaries}\label{sec2}
\def\theequation{2.\arabic{equation}}\makeatother
\setcounter{equation}{0}
\def\thetheorem{2.\arabic{theorem}}\makeatother
\setcounter{theorem}{0}

As a preparation for the proof of the global existence, we state the following lemmas which will be used often in the next sections.

The following Sobolev inequality is ensured in \cite{lieb202}.
\begin{lemma}\label{sobolev}
Let $n \ge 3$. Suppose $u \in H^1(\R^n)$. Then $u \in L^\frac{2n}{n-2}(\R^n)$ and the following holds:
\begin{align}\label{sobolev1}
S_n \|u\|_{L^\frac{2n}{n-2}(\R^n)}^2 \le \|\nabla u\|_{L^2(\R^n)}^2,~~S_n=\frac{n(n-2)}{4}2^{\frac{2}{n}}\pi^{1+\frac{1}{n}}\Gamma\left( \frac{n+1}{2}\right)^{-\frac{2}{n}}.
\end{align}
\end{lemma}
An immediate consequence of \er{sobolev1} together with the H\"{o}lder inequality implies
\begin{lemma}\label{gnssd}
Let $n \ge 3$. Assume $1 <\frac{b}{a}<\frac{2n}{a(n-2)}$ satisfying $\frac{b}{a}=\frac{2}{a}+\frac{2}{n}$, then for $w^{1/a} \in H^1(\R^n)$ and $w \in L^1(\R^n)$, it holds
\begin{align}\label{sd0}
\|w\|_{L^{b/a}(\R^n)}^{b/a} \le \frac{1}{S_n} \|\nabla w^{1/a}\|_{L^2(\R^n)}^2 \|w\|_{L^1(\R^n)}^{2/n}.
\end{align}
\end{lemma}
\noindent\textbf{ Proof.} We employ the H\"{o}lder inequality with $1 <\frac{b}{a}<\frac{2n}{a(n-2)}$ that
\begin{align}\label{sd1}
\|w^{1/a}\|_{L^b(\R^n)}^b \le \|w^{1/a}\|_{L^\frac{2n}{n-2}(\R^n)}^{2} \|w^{1/a}\|_{L^a(\R^n)}^{\frac{2}{n}a}
\end{align}
which follows since $\frac{b}{a}=\frac{2}{a}+\frac{2}{n}$. Therefore, the validity of \er{sd0} is checked by inserting \er{sobolev1} with $u=w^{1/a}$ into \er{sd1}. $\Box$

Furthermore, we show the general GNS inequality.

\begin{lemma}[Gagliardo-Nirenberg-Sobolev inequality] \label{gnsyoung}
Let $n\ge 3, 1<\frac{b}{a}<\frac{2n}{a(n-2)}$ and
$\frac{b}{a}<\frac{2}{a}+\frac{2}{n}$. Assume $w \in L_+^1(\R^n)$
and $ w^{1/a} \in H^1(\R^n)$ with $a>0$, then
\begin{align*}
    \|w\|_{L^{b/a}(\R^n)}^{b/a} \le \left(1-\frac{1}{\delta}\right) \delta^{-\frac{1}{\delta-1}} \left(S_n C_0\right)^{-\frac{1}{\delta-1}} \|w\|_{L^1(\R^n)}^{\gamma}+ C_0 \| \nabla w^{1/a} \|_{L^2(\R^n)}^2,
\end{align*}
where
$$
\delta=\frac{2\left( \frac{1}{a}- \frac{n-2}{2n}
\right)}{\frac{b}{a}-1},~~\\
\gamma=1+\frac{2(b-a)}{2a-(b-2)n},
$$
and $C_0$ is an arbitrarily positive constant.
\end{lemma}
\noindent\textbf{ Proof.}
Taking $u=w^{1/a}$ in \er{sobolev1} and applying the H\"{o}lder inequality
with $1<\frac{b}{a}<\frac{2n}{a(n-2)}$  yield
\begin{align*}
    \|w\|_{L^{b/a}(\R^n)} \le \|w\|_{L^1(\R^n)}^{1- \theta}
    \|w\|_{L^{\frac{2n}{a(n-2)}}(\R^n)}^{\theta}=\|w\|_{L^1(\R^n)}^{1-\theta}
    \|w^{1/a}\|_{L^{2n/(n-2)}(\R^n)}^{\theta a} \le S_n^{- \theta
    a/2}\|w\|_{L^1(\R^n)}^{1-\theta} \|\nabla w^{1/a} \|_{L^2(\R^n)}^{\theta a},
\end{align*}
whence follows
\begin{align} \label{grad2}
    \|w\|_{L^{b/a}(\R^n)}^{b/a} \le C(n) \|w\|_{L^1(\R^n)}^{\frac{b}{a}(1-\theta)
   } \| \nabla w^{1/a} \|_{L^2(\R^n)}^{b \theta},
\end{align}
where
\begin{align*}
\theta
=\frac{\frac{1}{a}-\frac{1}{b}}{\frac{1}{a}-\frac{n-2}{2n}},~~
C(n)=S_n^{-b \theta /2}.
\end{align*}
It is easy to verify that $ b \theta <2$ if
$\frac{b}{a}<\frac{2}{a}+\frac{2}{n}$. Therefore, by the Young inequality we
have that for arbitrary $B>0$
\begin{align*}
\|w\|_{L^{b/a}(\R^n)}^{b/a} & \le C(n) \|w\|_{L^1(\R^n)}^{\frac{b}{a}(1-\theta)
   } \| \nabla w^{1/a} \|_{L^2(\R^n)}^{b \theta} \\
   & =C(n) B \|w\|_{L^1(\R^n)}^{\frac{b}{a}(1-\theta)} B^{-1} \| \nabla w^{1/a} \|_{L^2(\R^n)}^{b \theta} \\
   & \le C(n)\frac{B^{\delta'}}{\delta'} \|w\|_{L^1(\R^n)}^{\frac{b}{a}(1-\theta)
\delta'} + C(n) \frac{B^{-\delta}}{\delta} \|\nabla w^{1/a} \|_{L^2(\R^n)}^{b
\theta \delta},
\end{align*}
where $\frac{1}{\delta'}+\frac{1}{\delta}=1$ and $\delta=\frac{2\left( \frac{1}{a}- \frac{n-2}{2n} \right)}{\frac{b}{a}-1}$ such
that
$$
b \theta \delta=2.
$$
Letting $C_0=C(n) \frac{B^{-\delta}}{\delta}$ we conclude the proof. $\Box$ \\

Moreover, we still need the following lemmas which have been proved in \cite{BL13}.
\begin{lemma}[\cite{BL13}]\label{BL14ode}
Assume $y(t) \ge 0$ is a $C^{1}$ function for $t>0$ satisfying
$$y'(t)\le \eta- \beta y(t)^p$$
for any $p>1, \eta > 0, \beta>0$, then
$y(t)$ has the following hyper-contractive property
\begin{align} \label{agle1}
    y(t) \le (\eta/\beta)^{1/p}+ \left(\frac{1}{\beta(p-1)t}\right)^{\frac{1}{p-1}},  \quad \mbox{ for any ~} t >
    0.
\end{align}
Furthermore, if $y(0)$ is bounded, then
\begin{align} \label{ytlessy0}
 y(t) \le \max \left( y(0),  (\eta/\beta)^{1/p} \right).
\end{align}
\end{lemma}
More generally, we have
\begin{lemma}[\cite{BL13}]\label{BL13ode}
Assume $f(t)\ge 0$ is a non-increasing function for $t>0. ~y(t) \ge 0$ is a $C^1$ function and satisfies
$$y'(t)\le f(t)- \beta y(t)^p$$
for any $p>1, \beta>0$, then for any $t_0>0$ one has
\begin{align}
y(t)\le \left( f(t_0)/\beta \right)^{1/p}+\left( \frac{1}{\beta (p-1)(t-t_0)} \right)^{\frac{1}{p-1}}, \quad \mbox{ for any ~} t >
    t_0.
\end{align}
\end{lemma}

\section{Global existence and asymptotic behavior in the case of $\alpha=1+2/n$}\label{sec3}
\def\theequation{3.\arabic{equation}}\makeatother
\setcounter{equation}{0}
\def\thetheorem{3.\arabic{theorem}}\makeatother
\setcounter{theorem}{0}

This section mainly focuses on the global existence of the classical solution in the case of $\alpha=1+2/n$. Starting from $M_0<C_{\alpha,C_\ast}$, where $C_{\alpha,C_\ast}$ is a universal constant depending on $\alpha$ and the bounded supremum of a functional (see Lemma \ref{cstarbdd}), Theorem \ref{critical} presents the existence of a global classical solution to \er{nkpp} with non-negative continuous initial data $U_0$ satisfying $U_0 \in L^1(\R^n) \cap L^\infty(\R^n)$. The long time behavior of the global solution is also obtained in Theorem \ref{longcritical}. Besides, with the purpose of relaxing the initial regularities for the global existence, the hyper-contractive estimates in Theorem \ref{hypercontractivity} deduce that when $U_0 \in L^1(\R^n)$, a global solution is bounded in $L^q(\R^n)$ for any $q>1$ and $t>0$.

\subsection{Global existence}

This subsection gives the existence of a time global classical solution to \er{nkpp}. We firstly establish a type of Gagliardo-Nirenberg-Sobolev inequality which will play an important role in deriving the a priori estimates of solutions to \er{nkpp}.
\begin{lemma}\label{cstarbdd}
Let $\alpha=1+2/n$. Assume $u \in L^1(\R^n)$ and $\nabla u \in L^2(\R^n)$, then
\begin{align}\label{cstar}
C_\ast:=\displaystyle \sup_{u \neq 0} \left\{ \frac{\|u\|_{L^{\alpha+1}(\R^n)}^{\alpha+1}}{\|u\|_{L^1(\R^n)}^{\alpha-1} \|\nabla u\|_{L^2(\R^n)}^2}, u \in L^1(\R^n), \nabla u \in L^2(\R^n), u \ge 0 \right\}<\infty.
\end{align}
\end{lemma}
\noindent\textbf{Proof.} Consider $u \in L^1(\R^n)$ and $\nabla u \in L^2(\R^n)$, we claim that $C_\ast$ is bounded from above by $\frac{1}{S_n}$ where $S_n$ is defined by \er{sobolev1}. Indeed, by Lemma \ref{sobolev} and the H\"{o}lder inequality with $\alpha=1+2/n$ we obtain
\begin{align*}
\|u\|_{L^{\alpha+1}(\R^n)}^{\alpha+1} \le \|u\|_{L^1(\R^n)}^{\alpha-1} \|u\|_{L^{\frac{2n}{n-2}}(\R^n)}^2 \le \frac{1}{S_n} \|u\|_{L^1(\R^n)}^{\alpha-1} \|\nabla u\|_{L^2(\R^n)}^2.
\end{align*}
Consequently, $C_\ast \le \frac{1}{S_n}$. $\Box$ \\

We next turn to the existence of the extremal of the GNS inequality which can be proved by similar arguments as for the variant HLS inequality in \cite{bj07}.
\begin{lemma}[The existence of $C_\ast$]\label{cstarexist}
Let $\alpha=1+2/n$. There exists a radially symmetric and non-increasing function $U \in L^1(\R^n)$ and $\nabla U \in L^2(\R^n)$ such that
\begin{align}
\|U\|_{L^{\alpha+1}(\R^n)}^{\alpha+1}=C_\ast \|U\|_{L^1(\R^n)}^{\alpha-1} \|\nabla U\|_{L^2(\R^n)}^2
\end{align}
with $\|U\|_{L^{\alpha+1}(\R^n)}=\|U\|_{L^1(\R^n)}=1$.
\end{lemma}
\noindent\textbf{ Proof.} We firstly define
\begin{align}
J(u):=\frac{\|u\|_{L^{\alpha+1}(\R^n)}^{\alpha+1}}{\|u\|_{L^1(\R^n)}^{\alpha-1} \|\nabla u\|_{L^2(\R^n)}^2},~~u \in L^1(\R^n),~\nabla u \in L^2(\R^n)
\end{align}
and a maximizing sequence $\{u_j\}$ with $u_j \in L^1(\R^n)$ and $\nabla u_j \in L^2(\R^n)$ such that
\begin{align}\label{maxcstar}
\displaystyle \lim_{j \to \infty} J(u_j)=C_\ast.
\end{align}
The proof can be divided into two steps. Firstly, we prove that the maximizing sequence $u_j$ can be assumed to be non-negative, radially symmetric and non-increasing with $\|u_j\|_{L^{\alpha+1}(\R^n)}=\|u_j\|_{L^1(\R^n)}=1$. The second step is devoted to guarantee that the supremum can be achieved and
\begin{align*}
\displaystyle \lim_{j \to \infty} u_j=U,~~~~C_\ast \|\nabla U\|_{L^2(\R^n)}^2=1
\end{align*}
with $\|U\|_{L^{\alpha+1}(\R^n)}=\|U\|_{L^1(\R^n)}=1$ and thus $J(U)=C_\ast.$

\noindent{\it\textbf{Step 1.}}(Radially symmetric, non-negative and non-increasing assumption) Actually, for any $\nabla u_j \in L^2(\R^n)$ one has
\begin{align*}
\int_{\R^n} \big|\nabla u_j\big|^2 dx=\int_{\R^n} \big|\nabla |u_j|\big|^2 dx
\end{align*}
which follows by \cite{GN98} such that
\begin{align}
J(|u_j|) = J(u_j)
\end{align}
provides that $|u_j|$ is also a maximizing sequence. Next we denote the scaling
\begin{align}
\overline{u}_j:=\lambda |u_j(\mu x) |
\end{align}
with
\begin{align*}
\left\{
  \begin{array}{ll}
    \mu=\|u_j\|_{L^1(\R^n)}^{\frac{\alpha+1}{n \alpha}}\|u_j\|_{L^{\alpha+1}(\R^n)}^{-\frac{\alpha+1}{n \alpha}}, \\
    \lambda=\mu^n \|u_j\|_{L^1(\R^n)}^{-1}.
  \end{array}
\right.
\end{align*}
Then we have
\begin{align*}
\|\overline{u}_j\|_{L^1(\R^n)}=\|\overline{u}_j\|_{L^{\alpha+1}(\R^n)}=1.
\end{align*}
A direct computation leads to
\begin{align}
J(\overline{u}_j)=J(|u_j|).
\end{align}
Finally denoting the symmetric non-increasing rearrangement of $\overline{u}_j$ by $u_j^\ast$, the Riesz's rearrangement inequalities
\cite[pp. 81]{lieb202} and \cite{BZ88} yield
\begin{align*}
&\|u_j^\ast\|_{L^1(\R^n)}=\|\overline{u}_j\|_{L^1(\R^n)}=1, \\
&\|u_j^\ast\|_{L^{\alpha+1}(\R^n)}=\|\overline{u}_j\|_{L^{\alpha+1}(\R^n)}=1, \\
&\|\nabla u_j^\ast\|_{L^2(\R^n)}\le \|\nabla \overline{u}_j\|_{L^2(\R^n)}
\end{align*}
whence follows
\begin{align}
J(u_j^\ast) \ge J(\overline{u}_j).
\end{align}
This entails that $u^\ast_j$ is also a maximizing sequence and thus we can assume that the maximizing sequence $u_j$ is a non-negative, radially symmetric and non-increasing function with $\|u_j\|_{L^1(\R^n)}=\|u_j\|_{L^{\alpha+1}(\R^n)}=1.$

\noindent{\it\textbf{Step 2.}}(Existence of the supremum) In this step, we will show that the maximizing sequence is convergent and the supremum of $J(u)$ can be achieved. In fact, by the assumption of $u_j$ together with the Sobolev inequality \er{sobolev1} we can give the following estimates
that for any $R>0$
\begin{align*}
\|u_j\|_{L^1(\R^n)} &= n \alpha_n \int_0^\infty u_j(r) r^{n-1}dr  \\
  & \ge n \alpha_n \int_0^R u_j(r) r^{n-1}dr \\
 & \ge \alpha_n u_j(R) R^n. \\
S_n^{-\frac{n}{n-2}} \|\nabla u_j\|_{L^2(\R^n)}^{\frac{2n}{n-2}} \ge \|u_j\|_{L^{\frac{2n}{n-2}}(\R^n)}^{\frac{2n}{n-2}}&= n \alpha_n \int_0^\infty u_j^{\frac{2n}{n-2}}(r) r^{n-1}dr  \\
& \ge n \alpha_n \int_0^R u_j^{\frac{2n}{n-2}}(r) r^{n-1}dr \\
& \ge \alpha_n u_j^{\frac{2n}{n-2}}(R) R^n.
\end{align*}
where $\alpha_n=|B(0,1)|$. So that
\begin{align*}
&u_j(R) \le \alpha_n^{-1} \|u_j\|_{L^1(\R^n)} R^{-n}, \\
&u_j(R) \le \alpha_n^{-\frac{n-2}{2n}} S_n^{-\frac{1}{2}} \|\nabla u_j\|_{L^2(\R^n)} R^{-\frac{n-2}{2}}.
\end{align*}
Therefore we have
\begin{align} \label{star00}
u_j(R) \le G(R):=C_0 \inf \left\{ R^{-n}, R^{-\frac{n-2}{2}} \right\}~~~~\mbox{for~any}~~R>0.
\end{align}
Now using the non-increasing of $u_j$ and their boundedness in $(R,\infty),$ the Helly's selection principle \cite[pp. 89]{lieb202} deduces that there are a sub-sequence of $u_j$ (not relabeled) such that
\begin{align}\label{star000}
u_j \to U~~~~\mbox{pointwisely},
\end{align}
where $U$ is a non-negative and non-increasing function. Besides, the fact $1<\alpha+1<\frac{2n}{n-2}$ yields
\begin{align*}
\|G(|x|)\|_{L^{\alpha+1}(\R^n)}^{\alpha+1} &=n \alpha_n \int_0^\infty G(r)^{\alpha+1} r^{n-1} dr \\
&=C_1 \left( \int_0^1 \frac{1}{r^{\frac{(n-2)(\alpha+1)}{2}}} r^{n-1} dr+\int_1^\infty \frac{1}{r^{n(\alpha+1)}} r^{n-1} dr  \right) \\
&< \infty.
\end{align*}
Together with \er{star00} and \er{star000} produces
\begin{align}
u_j \to U \mbox{~~in~~} L^{\alpha+1}(\R^n),~~~~\mbox{as~~} j \to \infty
\end{align}
by dominated convergence theorem. Thus one has
\begin{align}\label{Ualpha}
\displaystyle \lim_{j \to \infty} \|u_j\|_{L^{\alpha+1}(\R^n)} =\|U\|_{L^{\alpha+1}(\R^n)}=1.
\end{align}
In addition, \er{star000} and Fatou's lemma ensure
\begin{align*}
\|U\|_{L^1(\R^n)} &\le \displaystyle \liminf_{j \to \infty} \|u_j\|_{L^1(\R^n)}=1, \\
\|\nabla U\|_{L^2(\R^n)} &\le \displaystyle \liminf_{j \to \infty} \|\nabla u_j\|_{L^2(\R^n)}.
\end{align*}
Hence using \er{maxcstar} and \er{Ualpha} we conclude that
\begin{align}
C_\ast=\displaystyle \lim_{j \to \infty} J(u_j)=\displaystyle \lim_{j \to \infty} \frac{\|u_j\|_{L^{\alpha+1}(\R^n)}^{\alpha+1}}{ \|u_j\|_{L^1(\R^n)}^{\alpha-1} \|\nabla u_j\|_{L^2(\R^n)}^2} \le \frac{\|U\|_{L^{\alpha+1}(\R^n)}^{\alpha+1}}{ \|U\|_{L^1(\R^n)}^{\alpha-1} \|\nabla U\|_{L^2(\R^n)}^2}=J(U)\le C_\ast.
\end{align}
Thus $U$ is non-negative and non-increasing satisfying $J(U)=C_\ast$ with $\|U\|_{L^1(\R^n)}=\|U\|_{L^{\alpha+1}(\R^n)}=1.$ $\Box$ \\

Now we are in a position to state the main result:

\begin{theorem}[Time global existence of $\alpha=1+2/n$ case]\label{critical}
Suppose $\alpha=1+2/n$. Assume $U_0 \ge 0$ and
$$
U_0 \in C(\R^n)\cap L^\infty(\R^n),~~\|U_0\|_{L^1(\R^n)}<M_0
$$
with $M_0$ fulfilling
\begin{align}\label{M0}
\eta_0:=\left(\frac{\alpha-1}{C_\ast}\right)^{\frac{1}{\alpha}} \frac{\alpha}{\alpha-1}-M_0>0
\end{align}
where $C_\ast$ is defined by \er{cstar}, then problem \er{nkpp} possesses a unique and uniformly bounded non-negative global classical solution in $L^1 \cap L^\infty(\R^n)$ satisfying \er{property}. In addition,

\noindent\textbf{(i)} For any $1< k <\infty$, the global solution satisfies $m_0\le m(t) \le M_0$ and the following decay property that for any $t>0$
\begin{align}\label{decay12}
\|u(\cdot,t)\|_{L^k(\R^n)} \le C\left(\eta_0,m_0,k \right)~ t^{-\frac{k-1}{k(\alpha-1)}}
\end{align}
where $C\left(\eta_0,m_0,k\right)$ is a positive bounded constant depending on $\eta_0, m_0$ and $k$.

\noindent\textbf{(ii)} The global solution is uniformly bounded in time that
\begin{align}\label{uniformbdd}
\|u(\cdot,t)\|_{L^\infty(\R^n)} \le C\left( \|U_0\|_{L^1(\R^n)},\|U_0\|_{L^\infty(\R^n)} \right).
\end{align}
\end{theorem}

\noindent\textbf{Proof.} To simplify the presentation of the proof, we will therefore do the computations at a formal level for smooth solutions which behave well at infinity. The most important steps towards the global existence of classical solutions is the following steps 1-4. Then the standard parabolic theory with sufficient regularities allows us to state the result without further comment. Here the non-negativity of $u(x,t)$ is a consequence of the maximum principle since $U_0(x) \ge 0$.

We emphasize that $\|u\|_{L^2(\R^n)}$ will play a fundamental role in deriving the boundedness of $\|u\|_{L^k(\R^n)}$ for any $1<k\le \infty$. Firstly we give the reason.

\noindent{\it\textbf{Step 1.}}(The motivation of picking $k=2$ as a medium) We begin by multiplying \er{nkpp} with $k u^{k-1}(k>1)$ that
\begin{align}\label{guji1}
\frac{d}{dt} \int_{\R^n}  u^k dx + \frac{4(k-1)}{k} \int_{\R^n}
|\nabla u^{\frac{k}{2}} |^2 dx+ k \int_{\R^n} u dx \int_{\R^n} u^{k+\alpha-1}
dx = k M_0 \int_{\R^n} u^{k+\alpha-1} dx.
\end{align}
The crucial estimate is the following inequality. Using
\begin{align*}
w^{1/a}=u^{\frac{k}{2}},~~b=\frac{2(k+\alpha-1)}{k},~~a=\frac{n(\alpha-1)}{k}
\end{align*}
in Lemma \ref{gnssd} for any $k> \frac{(\alpha-1)(n-2)}{2}$ one has
\begin{align}\label{guji2}
\|u\|_{L^{k+\alpha-1}(\R^n)}^{k+\alpha-1} \le \frac{1}{S_n} \|\nabla u^{\frac{k}{2}}\|_{L^2(\R^n)}^2 \|u\|_{L^\frac{n(\alpha-1)}{2}(\R^n)}^{\alpha-1}.
\end{align}
Keeping in mind $\alpha=1+2/n$, substituting \er{guji2} into \er{guji1} we have
\begin{align}
\frac{d}{dt} \int_{\R^n}  u^k dx+k \int_{\R^n} u^{k+\alpha-1} dx\left( \frac{4(k-1)S_n}{k^2 \|u\|_{L^1(\R^n)}^{\alpha-1}}+\|u\|_{L^1(\R^n)}-M_0 \right)\le 0
\end{align}
for any $k \ge 1.$ Define
\begin{align}\label{guji4}
f(x)=\frac{4(k-1)S_n}{k^2 x^{\alpha-1}}+x-M_0,
\end{align}
after some computations we know that at
\begin{align}
x_0=\left( \frac{4(k-1)(\alpha-1)S_n}{k^2} \right)^\frac{1}{\alpha}
\end{align}
$f(x)$ attains its minimum
\begin{align}
f(x_0)=\left(\frac{4(k-1)(\alpha-1)S_n}{k^2}\right)^{\frac{1}{\alpha}} \frac{\alpha}{\alpha-1}-M_0.
\end{align}
In order to guarantee that $f(x)\ge f(x_0)$ is positive, we assume
\begin{align}\label{guji5}
M_0<\left(\frac{4(k-1)(\alpha-1)S_n}{k^2}\right)^{\frac{1}{\alpha}} \frac{\alpha}{\alpha-1}
\end{align}
such that
$$
f(x)\ge f(x_0)>0.
$$
Let's point out that for any $k>1$, $\frac{4(k-1)}{k^2}$ reaches its maximum at $k=2.$ Hence we firstly restrict the problem on $k=2$.

\noindent{\it\textbf{Step 2.}}(Decay estimates on $\|u\|_{L^2(\R^n)}$) Supposing $k=2$ in \er{guji1} one has that
\begin{align}\label{guji3}
\frac{d}{dt} \int_{\R^n}  u^2 dx + 2 \int_{\R^n}
|\nabla u |^2 dx+ 2 \int_{\R^n} u dx \int_{\R^n} u^{\alpha+1}
dx =2 M_0 \int_{\R^n} u^{\alpha+1} dx.
\end{align}
Thanks to Lemma \ref{cstarexist}, plugging
\begin{align}
\frac{\|u\|_{L^{\alpha+1}(\R^n)}^{\alpha+1}}{C_\ast \|u\|_{L^1(\R^n)}^{\alpha-1}} \le \|\nabla u\|_{L^2(\R^n)}^2
\end{align}
into \er{guji3} we proceed to derive
\begin{align}
\frac{d}{dt} \int_{\R^n}  u^2 dx+2 \left( \frac{1}{C_\ast \|u\|_{L^1(\R^n)}^{\alpha-1}}+\|u\|_{L^1(\R^n)}-M_0  \right) \|u\|_{L^{\alpha+1}(\R^n)}^{\alpha+1} \le 0.
\end{align}
In order to ensure that
$$
g\left( \|u\|_{L^1(\R^n)} \right)=\frac{1}{C_\ast \|u\|_{L^1(\R^n)}^{\alpha-1}}+\|u\|_{L^1(\R^n)}-M_0
$$
is positive, repeating the computations from \er{guji4} to \er{guji5} with $k=2$ we infer that
\begin{align}
M_0<\left(\frac{\alpha-1}{C_\ast}\right)^{\frac{1}{\alpha}} \frac{\alpha}{\alpha-1}
\end{align}
which results in
\begin{align*}
g\left( \|u\|_{L^1(\R^n)} \right) \ge g \left( \left( \frac{\alpha-1}{C_\ast} \right)^{\frac{1}{\alpha}} \right)=\left(\frac{\alpha-1}{C_\ast}\right)^{\frac{1}{\alpha}} \frac{\alpha}{\alpha-1}-M_0>0.
\end{align*}
So that
\begin{align}\label{guji6}
\frac{d}{dt} \int_{\R^n}  u^2 dx +2 \left( \left(\frac{\alpha-1}{C_\ast}\right)^{\frac{1}{\alpha}} \frac{\alpha}{\alpha-1}-M_0 \right) \|u\|_{L^{\alpha+1}(\R^n)}^{\alpha+1} \le 0.
\end{align}
Recall that the initial mass $m_0<M_0$, then \er{m00} implies
\begin{align}
\|u\|_{L^2(\R^n)}^{2 \alpha} \le \|u\|_{L^{\alpha+1}(\R^n)}^{\alpha+1} \|u\|_{L^1(\R^n)}^{\alpha-1} \le \|u\|_{L^{\alpha+1}(\R^n)}^{\alpha+1} M_0^{\alpha-1}
\end{align}
which follows by the H\"{o}lder inequality involving the fact $\alpha=1+2/n$. Then \er{guji6} becomes
\begin{align}
\frac{d}{dt} \int_{\R^n}  u^2 dx+ 2 \eta_0 \frac{\left(\int_{\R^n}  u^2 dx\right)^\alpha }{M_0^{\alpha-1}} \le 0
\end{align}
where $\eta_0=\left(\frac{\alpha-1}{C_\ast}\right)^{\frac{1}{\alpha}} \frac{\alpha}{\alpha-1}-M_0>0.$ Consequently it yields
\begin{align}\label{decay2}
\int_{\R^n} u^2 dx \le \left( \frac{1}{\left( \int_{\R^n} U_0^2 dx \right)^{1-\alpha}+(\alpha-1)C(\eta_0,M_0)~ t}  \right)^{\frac{1}{\alpha-1}}
\end{align}
which immediately leads to \er{decay12} as a consequence of $\|u\|_{L^k(\R^n)}^k \le \|u\|_{L^1(\R^n)}^{2-k} \|u\|_{L^2(\R^n)}^{2(k-1)}$ saturating the H\"{o}lder inequality with $1<k<2$.

\noindent{\it\textbf{Step 3.}}(Decay estimates on $\|u\|_{L^k(\R^n)}$ for $2<k<\infty$) Concerning $k>2$ and inserting
\begin{align}
w^{1/a}=u^{k/2}, ~~b=\frac{2(k+\alpha-1)}{k},~~ a=\frac{4}{k},~~ C_0=\frac{4(k-1)}{k^2 M_0}
\end{align}
into Lemma \ref{gnsyoung} for any $k>\frac{(\alpha-1)(n-2)}{2}=\frac{n-2}{n}$ we compute
\begin{align*}
\gamma=1+\frac{2(b-a)}{2a-n(b-2)}=k+\alpha-2 \mbox{~~for~~} \alpha=1+\frac{2}{n}
\end{align*}
and
\begin{align*}
\int_{\R^n} u^{k+\alpha-1} dx \le \frac{4(k-1)}{k^2 M_0} \| \nabla u^{k/2} \|_{L^2(\R^n)}^2 +C(k,M_0) \|u\|_{L^2(\R^n)}^{2 \gamma}
\end{align*}
for any $k \ge 1.$ Together with \er{guji1} we have
\begin{align}
\frac{d}{dt} \int_{\R^n}  u^k dx +k \int_{\R^n} u dx \int_{\R^n} u^{k+\alpha-1}
dx \le C(k,M_0) \left(\|u\|_{L^2(\R^n)}^2 \right)^{k+\alpha-2}.
\end{align}
Recalling $m_0 \le m(t) \le M_0$ with
\begin{align}\label{k1chazhi}
\left(\|u\|_{L^k(\R^n)}^{k}\right)^{1+\frac{\alpha-1}{k-1}} \le \|u\|_{L^{k+\alpha-1}(\R^n)}^{k+\alpha-1} \|u\|_{L^1(\R^n)}^{\frac{\alpha-1}{k-1}}\le \|u\|_{L^{k+\alpha-1}(\R^n)}^{k+\alpha-1} M_0^{\frac{\alpha-1}{k-1}}
\end{align}
allows us to get
\begin{align}\label{guji7}
\frac{d}{dt} \int_{\R^n}  u^k dx +C(m_0,M_0,k) \left(\int_{\R^n} u^k dx\right)^{1+\frac{\alpha-1}{k-1}} & \le C(k,M_0)~ \left(\|u\|_{L^2(\R^n)}^2\right)^{k+\alpha-2} \nonumber\\
& \le C(k,M_0,\eta_0) ~t^{-\frac{k+\alpha-2}{\alpha-1}}
\end{align}
where the last inequality follows by \er{decay2}. Setting
\begin{align*}
y(t)= \int_{\R^n}  u^k dx,~~ p=1+\frac{\alpha-1}{k-1},~~ f(t)=C(k,M_0,\eta_0)~ t^{-\frac{k+\alpha-2}{\alpha-1}}
\end{align*}
in Lemma \ref{BL13ode} one obtains that for $t_0=t/2$ with any $t>0$
\begin{align}
\|u\|_{L^k(\R^n)}^k \le C\left(k,m_0,\eta_0 \right) t^{-\frac{k-1}{\alpha-1}},~~\mbox{for any}~~2< k<\infty.
\end{align}
The decay rate in time is consistent with \er{decay2} and thus we conclude \er{decay12} for any $2<k<\infty.$ This allows us to go further to estimate $\|u\|_{L^\infty(\R^n)}.$

\noindent{\it\textbf{Step 4.}}(Uniformly boundedness) The goal of this step is to derive the uniformly boundedness in time of the solution. We denote
\begin{align}
q_m=2^m+\alpha,~~m \ge 0.
\end{align}
Multiplying $q_m u^{q_m-1}$ to \er{nkpp} one obtains
\begin{align}\label{100}
&\frac{d}{dt} \int_{\R^n}  u^{q_m} dx + \frac{4(q_m-1)}{q_m} \int_{\R^n}
|\nabla u^{\frac{q_m}{2}} |^2 dx+ q_m \int_{\R^n} u dx \int_{\R^n} u^{q_m+\alpha-1}
dx \nonumber\\
=  &q_m M_0\int_{\R^n} u^{q_m+\alpha-1} dx.
\end{align}
Armed with Lemma \ref{gnsyoung} and letting
\begin{align*}
w^{\frac{1}{a}}=u^{\frac{q_m}{2}},~~b=\frac{2(q_m+\alpha-1)}{q_m},~~a=\frac{2q_{m-1}}{q_m},~~C_0=\frac{1}{q_m M_0 }
\end{align*}
we estimate
\begin{align}
q_m M_0\int_{\R^n} u^{q_m+\alpha-1} dx \le \left(\delta_1 S_n\right)^{-\frac{1}{\delta_1-1}} M_0^{\frac{\delta_1}{\delta_1-1}} q_m^{\frac{\delta_1}{\delta_1-1}} \left( \int_{\R^n} u^{q_{m-1}} dx \right)^{\gamma_1}+\|\nabla u^{\frac{q_m}{2}} \|_{L^2(\R^n)}^2
\end{align}
where
\begin{align*}
&\delta_1=\frac{q_m-\frac{n-2}{n}q_{m-1}}{q_m+\alpha-1-q_{m-1}}=\alpha, \\
&\gamma_1=1+\frac{2(b-a)}{2a-(b-2)n}=2
\end{align*}
which were ensured by $\alpha=1+2/n.$ Substituting the above into \er{100} and noticing that
$$
\frac{4(q_m-1)}{q_m} >2
$$
follow that
\begin{align}\label{200}
\frac{d}{dt} \int_{\R^n}  u^{q_m} dx + \int_{\R^n}
|\nabla u^{\frac{q_m}{2}} |^2 dx+ q_m \int_{\R^n} u dx \int_{\R^n} u^{q_m+\alpha-1}
dx \le C(\alpha,M_0) q_m^{\frac{\alpha}{\alpha-1}} \left( \int_{\R^n} u^{q_{m-1}} dx \right)^2.
\end{align}
Again applying Lemma \ref{gnsyoung} with
\begin{align*}
w^{\frac{1}{a}}=u^{\frac{q_m}{2}},~~b=2,~~a=\frac{2q_{m-1}}{q_m},~~C_0=1
\end{align*}
and using the Young inequality we have
\begin{align}\label{300}
\int_{\R^n} u^{q_m}dx &\le C(n) \left( \int_{\R^n}u^{q_{m-1}} dx \right)^{\gamma_2}+\|\nabla u^{\frac{q_m}{2}} \|_{L^2(\R^n)}^2 \nonumber \\
& \le \alpha \int_{\R^n} u dx \int_{\R^n} u^{q_m+\alpha-1}dx+C(\alpha)+\|\nabla u^{\frac{q_m}{2}} \|_{L^2(\R^n)}^2
\end{align}
where we have used
\begin{align*}
\gamma_2=1+\frac{q_m-q_{m-1}}{q_{m-1}}<2,~~q_{m-1}=\frac{q_m+\alpha}{2}.
\end{align*}
Plugging \er{300} into \er{200} one has
\begin{align*}
\frac{d}{dt} \int_{\R^n}  u^{q_m} dx + \int_{\R^n}  u^{q_m} dx & \le C_1(\alpha,M_0) q_m^{\frac{\alpha}{\alpha-1}} \left( \int_{\R^n} u^{q_{m-1}} dx \right)^2+C_2(\alpha) \\
& \le \max\left\{C_1(\alpha,M_0),C_2(\alpha)  \right\} 2^{\frac{\alpha}{\alpha-1}m} \max\left\{ 1, \left( \int_{\R^n} u^{q_{m-1}} dx \right)^2 \right\} \\
& = C(\alpha,M_0) 2^{\frac{\alpha}{\alpha-1}m} \max\left\{ 1, \left( \int_{\R^n} u^{q_{m-1}} dx \right)^2 \right\}
\end{align*}
Denote $y_m(t)=\int_{\R^n} u^{q_m} dx,$ then it satisfies
\begin{align*}
\left( e^t y_m(t) \right)' \le C(\alpha,M_0) 2^{\frac{\alpha}{\alpha-1}m} \max\left\{ 1, y_{m-1}^2(t) \right\}~e^t
\end{align*}
whence we find
\begin{align}\label{400}
y_m(t) & \le \left(1-e^{-t}\right) C(\alpha,M_0) 2^{\frac{\alpha}{\alpha-1}m} \max\left\{1,\displaystyle \sup_{ t \ge 0} y_{m-1}^2(t)  \right\}+e^{-t} y_m(0) \nonumber \\
& \le \max\left\{C(\alpha,M_0),1\right\}~2^{\frac{\alpha}{\alpha-1}m} \max\left\{1,\displaystyle \sup_{ t \ge 0} y_{m-1}^2(t),y_m(0)\right\} \nonumber \\
& =\overline{C}(\alpha,M_0) 2^{\frac{\alpha}{\alpha-1}m} \max\left\{1,\displaystyle \sup_{ t \ge 0} y_{m-1}^2(t),y_m(0)\right\}.
\end{align}
Define $\overline{K}=\max\big\{1,\|U_0\|_{L^1(\R^n)},\|U_0\|_{L^\infty(\R^n)}\big\}$, we have the following inequality for the initial data
\begin{align*}
y_m(0)=\int_{\R^n} U_0^{q_m} dx \le \overline{K}^{q_m} \le \overline{K}^{(\alpha+1) 2^m}=K_0^{2^m}.
\end{align*}
Thus from \er{400} one obtains
\begin{align*}
y_m(t) \le \overline{C}(\alpha,M_0) 2^{\frac{\alpha}{\alpha-1}m} \max\left\{K_0^{2^m}, \displaystyle \sup_{ t \ge 0} y_{m-1}^2(t) \right\}.
\end{align*}
After some iterative steps we can estimate
\begin{align*}
y_m(t) & \le \left(\overline{C}(\alpha,M_0)\right)^{1+2+\cdots 2^{m-1}} 2^{\frac{\alpha}{\alpha-1}\left(m+2(m-1)+2^2(m-2)+\cdots 2^{m-1}\right)} \max\left\{K_0^{2^m}, \displaystyle \sup_{ t \ge 0} y_{0}^{2^m}(t) \right\} \\
&= \left(\overline{C}(\alpha,M_0)\right)^{2^m-1} 2^{\frac{\alpha}{\alpha-1}\left( 2^{m+1}-m-2\right)} \max\left\{K_0^{2^m}, \displaystyle \sup_{ t \ge 0} y_{0}^{2^m}(t) \right\}.
\end{align*}
It's equivalent to
\begin{align}\label{500}
y_m(t)\le \left(\overline{C}(\alpha,M_0)\right)^{2^m-1} 2^{\frac{\alpha}{\alpha-1}\left( 2^{m+1}-m-2\right)} \max\left\{K_0^{2^m}, \displaystyle \sup_{ t \ge 0} \left( \int_{\R^n} u^{q_0} dx \right)^{2^m} \right\}.
\end{align}
Taking the power $\frac{1}{q_m}$ to both sides of \er{500}, then the uniformly boundedness of the global solution is obtained by passing to the limit $m \to \infty$
\begin{align}\label{uniform00}
\|u(\cdot,t)\|_{L^\infty(\R^n)} \le \overline{C}(\alpha,M_0) 2^{\frac{2\alpha}{\alpha-1}} \max\left\{K_0, \displaystyle \sup_{ t \ge 0} \int_{\R^n} u^{q_0}(t) dx  \right\}.
\end{align}
Now we turn to estimate $\int_{\R^n} u^{q_0} dx=\int_{\R^n} u^{\alpha+1} dx.$ Back to \er{guji7}, recalling \er{decay2} one has that for any $k>2$
\begin{align*}
\frac{d}{dt} \int_{\R^n}  u^k dx +C(m_0,M_0,k) \left(\int_{\R^n} u^k dx\right)^{1+\frac{\alpha-1}{k-1}} & \le C(k,M_0) \|u\|_{L^2(\R^n)}^{2 (k+\alpha-2)} \\
& \le C(k,M_0) \|U_0\|_{L^2(\R^n)}^{2 (k+\alpha-2)}.
\end{align*}
Letting $k=\alpha+1,$ plugging
$$
\beta=C(m_0,M_0,k),~~\eta=C(k,M_0) \|U_0\|_{L^2(\R^n)}^{2 (k+\alpha-2)},~~p=1+\frac{\alpha-1}{k-1}
$$
into Lemma \ref{BL14ode} gives
\begin{align}
\int_{\R^n} u^{\alpha+1} dx \le \max\left\{ \int_{\R^n} U_0^{\alpha+1} dx,  C(m_0,M_0,\alpha) \|U_0\|_{L^2(\R^n)}^{2\alpha} \right\}.
\end{align}
Hence \er{uniform00} becomes
\begin{align}
\|u(\cdot,t)\|_{L^\infty(\R^n)} \le C\left( \|U_0\|_{L^1(\R^n)},\|U_0\|_{L^\infty(\R^n)} \right).
\end{align}

\noindent{\it\textbf{Step 5.}}(Global existence) We now have necessary a priori estimates for the existence of global classical solutions. We know that $u$ is uniformly bounded for any $t \ge 0$ and the reaction term $u^\alpha\left( M_0 -\int_{\R^n} u dx\right)$ is bounded from below and above. Hence the global existence of classical solutions is followed by the standard parabolic theory for semilinear equations. In the end, the uniqueness can be obtained from the comparison principle, since $u^\alpha \left( M_0-\int_{\R^n} u dx \right)$ is bounded from below and above. This completes the proof of the global existence and uniqueness of the classical solution. $\Box$

\subsection{Asymptotic behavior of the global solution} \label{subsec2}
We study the large time behavior of the global solution via the corresponding integral equation to \er{nkpp}
\begin{align}\label{kernel}
u(x,t)=e^{t\Delta} U_0 (x)+\int_0^t e^{(t-s)\Delta} u(s)^\alpha \left( M_0-\int_{\R^n} u(y,s)dy \right)  ds.
\end{align}
Recall that
\begin{align}
e^{t\Delta}U_0(x)=\int_{\R^n} G_t(x-y) U_0(y) dy,~~G_t(x)=\frac{1}{(4\pi t)^{\frac{n}{2}}} e^{-\frac{|x|^2}{4t}}.
\end{align}
For future reference we collect some well known facts \cite{weis81} about the semigroup $e^{t\Delta}$.

\begin{proposition}\label{prop1}
\begin{itemize}
  \item[(a)] $\|G_t\|_{L^1(\R^n)}=1$ for all $t>0.$
  \item[(b)] If $v \ge 0,$ then $e^{t\Delta} v \ge 0$ and $\|e^{t \Delta} v\|_{L^1(\R^n)}=\|v\|_{L^1(\R^n)}$.
  \item[(c)] If $1\le p\le \infty,$ then $\|e^{t\Delta} v\|_{L^p(\R^n)}\le \|v\|_{L^p(\R^n)}$ for all $t>0.$
  \item[(d)] If $1\le p\le q\le \infty,$ then $\|e^{t\Delta}v\|_{L^q(\R^n)}\le \frac{1}{(4\pi t)^{\frac{n}{2} \left( \frac{1}{p}-\frac{1}{q} \right)}} \|v\|_{L^p(\R^n)}$ for all $t>0.$
\end{itemize}
\end{proposition}

We are now state the long time behavior of the global solution.
\begin{theorem}[Asymptotic behavior of the global solution for $\alpha=1+2/n$]\label{longcritical}
Suppose $\alpha=1+2/n$. $u(x,t)$ is the global solution to \er{nkpp}, then
\begin{align}
&\|u(\cdot,t)-G_t\ast U_0(x) \|_{L^\infty(\R^n)} \nonumber \\
\le& B\left( 1-\frac{n(r-1)}{2p}, 1-\frac{n}{2p} \right) C\left( \|U_0\|_{L^1(\R^n)},\|U_0\|_{L^\infty(\R^n)} \right) \left( M_0-m_0 \right) ~t^{-\left( \frac{nr}{2p}-1\right)}
\end{align}
for any $1<\frac{2p}{n}<r<\frac{2p}{n}+1.$
\end{theorem}
\noindent\textbf{ Proof.} From the integral equation \er{kernel} with the initial value $U_0(x) \ge 0$ and by Proposition \ref{prop1} it follows that for any $1 < p<q \le \infty$
\begin{align}
\|u(\cdot,t)-G_t\ast U_0(x) \|_{L^q(\R^n)} & \le \int_0^t \left\| \int_0^t e^{(t-s)\Delta} u(s)^\alpha \left( M_0-\int_{\R^n} u(y,s)dy \right) \right\|_{L^q(\R^n)}ds \nonumber\\
& \le \int_0^t [4\pi(t-s)]^{-\frac{n}{2}\left( \frac{1}{p}-\frac{1}{q} \right)} \|u(s)\|_{L^{p \alpha}(\R^n)}^\alpha \left( M_0-\int_{\R^n} u(y,s)dy \right) ds \nonumber\\
& \le \int_0^t [4\pi(t-s)]^{-\frac{n}{2}\left( \frac{1}{p}-\frac{1}{q} \right)} \|u(s)\|_{L^{p \alpha}(\R^n)}^\alpha \left( M_0-m_0 \right) ds. \label{guji9}
\end{align}
The last inequality follows since $M_0-m(t) \le M_0-m_0$ by \er{m00}. Now we choose $1<r<p\alpha$ such that
\begin{align*}
\|u(s)\|_{L^{p \alpha}(\R^n)}^\alpha &=\left(\int_{\R^n} u^{p\alpha-r} u^r dx\right)^{\frac{1}{p}} \\
& \le \|u\|_{L^\infty(\R^n)}^{\frac{p\alpha-r}{p}} s^{-\frac{r-1}{(\alpha-1)p}}
\end{align*}
in light of \er{decay12}. Then \er{uniformbdd} and the choice of $r$ guarantee that
\begin{align*}
& \|u(\cdot,t)-G_t\ast U_0(x) \|_{L^q(\R^n)} \le C\left( \|U_0\|_{L^1(\R^n)},\|U_0\|_{L^\infty(\R^n)} \right) \left( M_0-m_0 \right)\int_0^t (t-s)^{-\frac{n}{2}\left( \frac{1}{p}-\frac{1}{q} \right)} s^{-\frac{r-1}{(\alpha-1)p}}ds \nonumber\\
& = C\left( \|U_0\|_{L^1(\R^n)},\|U_0\|_{L^\infty(\R^n)} \right) \left( M_0-m_0 \right) t^{1-\frac{n}{2}\left( \frac{1}{p}-\frac{1}{q} \right)-\frac{r-1}{(\alpha-1)p}} B\left( 1- \frac{r-1}{(\alpha-1)p}, 1-\frac{n}{2}\left( \frac{1}{p}-\frac{1}{q} \right) \right).
\end{align*}
Let $q=\infty,$ taking $\alpha=1+2/n$ into account together with
\begin{align}\label{191012}
 \frac{nr}{2p}>1,~~\frac{n}{2p}(r-1)<1,~~\frac{n}{2p}<1
\end{align}
it suffices to show that
\begin{align*}
\|u(\cdot,t)-G_t\ast U_0(x) \|_{L^\infty(\R^n)} \le B\left( 1-\frac{n(r-1)}{2p}, 1-\frac{n}{2p} \right) C\left( \|U_0\|_{L^1(\R^n)},\|U_0\|_{L^\infty(\R^n)} \right) \left( M_0-m_0 \right) ~t^{-\left( \frac{nr}{2p}-1\right)}
\end{align*}
for any $1<\frac{2p}{n}<r<\frac{2p}{n}+1$. Here $B\left( 1-\frac{n(r-1)}{2p}, 1-\frac{n}{2p} \right)$ is bounded from above with the choice \er{191012}.   This ends up with the convergence to the solution of the heat equation as desired. $\Box$

\subsection{Hyper-contractivity of a global solution}\label{subsec3}

By Theorem \ref{critical} we know that there is a global classical solution which is uniformly bounded in time when $U_0 \in L^1 \cap L^\infty(\R^n)$. However, we wouldn't like to make such a restrictive regularity on the initial data. In fact, to an initial data $U_0(x)$ which is originally in $L^1(\R^n)$ but not in $L^\infty(\R^n)$, the following theorem associates a global solution which is in $L^k(\R^n)$ for any $k>1$ and thus gives the hyper-contractivity (defined in \cite[pp. 24]{bdp06}) of the global solution.

\begin{theorem}[Hyper-contractivity of a global solution] \label{hypercontractivity}
Suppose $\alpha=1+2/n$. Assume the non-negative and non-trivial initial data satisfies
\begin{align}\label{M00}
\|U_0\|_{L^1(\R^n)}<M_0,~~\eta_0:=\left(\frac{\alpha-1}{C_\ast}\right)^{\frac{1}{\alpha}} \frac{\alpha}{\alpha-1}-M_0>0
\end{align}
where $C_\ast$ is defined by \er{cstar}. It follows that $u(\cdot,t) \in L^1(\R^n)$ for any $t \in [0,\infty)$ and $u(\cdot,t) \in L^\alpha(\R^n)$ for any $t \in (0,\infty)$ which is a non-negative solution of the integral equation \er{property}. Moreover,

\noindent\textbf{(i)} For $1<k <\infty$, the following hyper-contractive estimates of the global solution hold true that for any $t>0$
\begin{align}\label{decay120}
\|u(\cdot,t)\|_{L^k(\R^n)} \le C(\eta_0,m_0,k)~ t^{-\frac{k-1}{k(\alpha-1)}}.
\end{align}

\noindent\textbf{(ii)} Furthermore, the global solution fulfilling ultra-contractivity.
\begin{align}
\|u\|_{L^\infty(\R^n)} \le C\left(\eta_0,m_0\right) ~t^{-\frac{\alpha}{\alpha-1}+\frac{n}{2}},~~\mbox{for~any}~~0<t\le 1.
\end{align}
and
\begin{align}
\|u\|_{L^\infty(\R^n)} \le C(\eta_0,m_0)~ t^{-\frac{\alpha}{\alpha-1}},~~\mbox{for~any}~~1<t<\infty.
\end{align}
\end{theorem}
\noindent\textbf{Proof.} The global existence of a global solution has been proved in Theorem \ref{critical}. Now we will give the hyper-contractivity in $L^k(\R^n)$ for any $ 1<k \le \infty$.

Repeating the process from step 1 to step 3 in the proof of Theorem \ref{critical} allows us to get that for any $t>0$ and $1<k<\infty$
\begin{align}\label{0918}
\|u(\cdot,t)\|_{L^k(\R^n)}^k \le C(\eta_0,m_0,k)~t^{-\frac{k-1}{\alpha-1}}.
\end{align}
Similar to the computations from \er{100} to \er{200} we also obtain
\begin{align}\label{0913}
\frac{d}{dt} \int_{\R^n}  u^{q_m} dx + \int_{\R^n}
|\nabla u^{\frac{q_m}{2}} |^2 dx+ q_m \int_{\R^n} u dx \int_{\R^n} u^{q_m+\alpha-1}
dx \le C(\alpha,M_0) q_m^{\frac{\alpha}{\alpha-1}} \left( \int_{\R^n} u^{q_{m-1}} dx \right)^2
\end{align}
where $q_m=2^m+\alpha$ with $m \ge 0.$ On the other hand, letting
\begin{align*}
w^{\frac{1}{a}}=u^{\frac{q_m}{2}},~~b=2,~~a=\frac{2q_{m-1}}{q_m}
\end{align*}
in Lemma \ref{gnsyoung} we get that
\begin{align*}
\left( \|u\|_{L^{q_m}(\R^n)}^{q_m} \right)^{1+ \frac{2q_{m-1}/n}{q_m-q_{m-1}}} \le S_n^{-1} \left( \|u\|_{L^{q_{m-1}}(\R^n)}^{q_{m-1}} \right)^{\frac{2q_m/n}{q_m-q_{m-1}}} \|\nabla u^{\frac{q_m}{2}} \|_{L^2(\R^n)}^2.
\end{align*}
Substituting it into \er{0913} follows that for any $t>0$
\begin{align}
\frac{d}{dt}\int_{\R^n} u^{q_m} dx & \le -\frac{S_n}{ \left( \|u\|_{L^{q_{m-1}}(\R^n)}^{q_{m-1}} \right)^{\frac{2q_m/n}{q_m-q_{m-1}}} } \left( \|u\|_{L^{q_m}(\R^n)}^{q_m} \right)^{1+ \frac{2q_{m-1}/n}{q_m-q_{m-1}}}+C(\alpha) 2^{\frac{\alpha}{\alpha-1}m} \left( \int_{\R^n} u^{q_{m-1}} dx \right)^2 \nonumber\\
& \le -\frac{S_n}{ \displaystyle \sup_{0<t<\infty} \left( \int_{\R^n} u^{q_{m-1}}dx \right)^{\frac{2q_m/n}{q_m-q_{m-1}}} } \left( \|u\|_{L^{q_m}(\R^n)}^{q_m} \right)^{1+ \frac{2q_{m-1}/n}{q_m-q_{m-1}}}+C(\alpha) 2^{\frac{\alpha m}{\alpha-1}} \displaystyle \sup_{0<t<\infty} \left( \int_{\R^n} u^{q_{m-1}} dx \right)^2.
\end{align}
Thus plugging
\begin{align*}
&p=1+ \frac{2q_{m-1}/n}{q_m-q_{m-1}},~~\beta=\frac{S_n}{ \displaystyle \sup_{0<t<\infty} \left( \int_{\R^n} u^{q_{m-1}} dx \right)^{\frac{2q_m/n}{q_m-q_{m-1}}} } ,~~\\
&\eta=C(\alpha) 2^{\frac{\alpha m}{\alpha-1}} \displaystyle \sup_{0<t<\infty} \left( \int_{\R^n} u^{q_{m-1}} dx \right)^2
\end{align*}
into Lemma \ref{BL14ode} and letting $y_m(t)=\int_{\R^n} u^{q_m} dx$ one has
\begin{align*}
y_m(t) & \le \left( \frac{C(\alpha)}{S_n} \right)^{\frac{1}{p}} 2^{\frac{\alpha m}{ \alpha-1} \frac{1}{p}} \displaystyle \sup_{0<t<\infty} \left( \int_{\R^n} u^{q_{m-1}} dx \right)^{ \frac{1+ \frac{2q_{m}/n}{q_m-q_{m-1}}}{1+ \frac{2q_{m-1}/n}{q_m-q_{m-1}}} } +C(n)^{\frac{q_m-q_{m-1}}{2q_{m-1}/n}} \frac{1}{t^{\frac{q_m-q_{m-1}}{2q_{m-1}/n}}} \displaystyle \sup_{0<t<\infty} \left( \int_{\R^n} u^{q_{m-1}} dx \right)^{\frac{q_m}{q_{m-1}}} \\
& \le \max\left\{ 1, \frac{C(\alpha)}{S_n},C(n)  \right\} 2^{\frac{\alpha}{\alpha-1}m} \left( \displaystyle \sup_{0<t<\infty} \left( \int_{\R^n} u^{q_{m-1}} dx \right)^A+ \frac{1}{t^{\frac{q_m-q_{m-1}}{2q_{m-1}/n}}} \displaystyle \sup_{0<t<\infty} \left( \int_{\R^n} u^{q_{m-1}} dx \right)^B \right)
\end{align*}
where
\begin{align*}
A=\frac{2+ \frac{2q_{m}/n}{q_m-q_{m-1}}}{1+ \frac{2q_{m-1}/n}{q_m-q_{m-1}}}<2,~~B=\frac{q_m}{q_{m-1}}<2.
\end{align*}
Here we have used $p>1$. Denote $C_n=\max\left\{ 1, \frac{C(\alpha)}{S_n},C(n)  \right\} $, it allows us to go further that
\begin{align}\label{0914}
y_m(t)\le C_n 2^{\frac{\alpha m}{\alpha-1}} \left\{ \max\left( 1, \displaystyle \sup_{0<t<\infty} \left( \int_{\R^n} u^{q_{m-1}} dx \right)^2 \right) +   \frac{1}{t^{\frac{q_m-q_{m-1}}{2q_{m-1}/n}}} \max\left( 1, \displaystyle \sup_{0<t<\infty} \left( \int_{\R^n} u^{q_{m-1}} dx \right)^2 \right)  \right\}.
\end{align}

If $0<t \le 1,$ then $\frac{q_m-q_{m-1}}{2q_{m-1}/n}<\frac{n}{2}$ gives rise to
\begin{align*}
y_m(t) \le \frac{2C_n}{t^{n/2}} 2^{\frac{\alpha}{\alpha-1}m}  \max\left( 1, \displaystyle \sup_{0<t<\infty} y_{m-1}^2(t) \right).
\end{align*}
Some iterative procedures deduce that for any $0<t<1$
\begin{align*}
y_m(t) \le \left(\frac{2C_n}{t^{n/2}}\right)^{2^m-1} 2^{\frac{\alpha}{\alpha-1}\left( 2^{m+1}-m-2 \right)} \max\left( 1, \displaystyle \sup_{0<t<\infty} y_0^{2^m}(t) \right).
\end{align*}
Taking the power $\frac{1}{q_m}$ to both sides we conclude that for any $0<t<1$
\begin{align}
\|u(\cdot,t)\|_{L^\infty(\R^n)} &\le \frac{2C_n}{t^{n/2}} 2^{\frac{2 \alpha}{\alpha-1}} \max\left( \displaystyle \sup_{0<t<\infty} y_0(t),1 \right) \nonumber \\
&=\frac{2C_n}{t^{n/2}} 2^{\frac{2 \alpha}{\alpha-1}} \max\left( \displaystyle \sup_{0<t<\infty} \|u(\cdot,t)\|_{L^{\alpha+1}(\R^n)}^{\alpha+1},1 \right)
\nonumber \\
& \le 2C_n 2^{\frac{2 \alpha}{\alpha-1}}  \max\{ C(\eta_0,m_0,\alpha),1\} ~t^{-\frac{\alpha}{\alpha-1}+\frac{n}{2}}
\end{align}
where we have applied \er{0918}. Similarly, if $1<t<\infty,$ then \er{0914} becomes
\begin{align*}
y_m(t) \le 2C_n 2^{\frac{\alpha}{\alpha-1}m}  \max\left( 1, \displaystyle \sup_{0<t<\infty} y_{m-1}^2(t) \right)
\end{align*}
whence follows
\begin{align*}
\|u(\cdot,t)\|_{L^\infty(\R^n)} & \le 2 C_n 2^{2\frac{\alpha}{\alpha-1}} \max\left( \displaystyle \sup_{0<t<\infty} \|u(\cdot,t)\|_{L^{\alpha+1}(\R^n)}^{\alpha+1},1 \right)
& \le 2 C_n 2^{2\frac{\alpha}{\alpha-1}}  \max\{ C(\eta_0,m_0,\alpha),1\}~ t^{-\frac{\alpha}{\alpha-1}}.
\end{align*}
We now have necessary a priori estimates for the global existence and uniqueness of the solution to \er{nkpp} for any $t>0$, which is followed by the standard parabolic theory. Thus completes the proof. $\Box$

\section{Global existence for $1<\alpha<1+2/n$}\label{sec4}
\def\theequation{4.\arabic{equation}}\makeatother
\setcounter{equation}{0}
\def\thetheorem{4.\arabic{theorem}}\makeatother
\setcounter{theorem}{0}

In this section we shall construct a global solution to \er{nkpp} for $1<\alpha<1+2/n.$ The result shows that the solution exists globally without any restriction on $M_0$ which is in sharp contrast to the case $\alpha=1+2/n.$ Precisely, we have:

\begin{theorem}[Time global existence of $1<\alpha<1+2/n$ case]\label{subcritical}
Suppose $1<\alpha<1+2/n$. Assume $U_0(x)$ is a non-negative bounded continuous function and $m_0=\int_{\Omega} U_0(x) dx<M_0$, then
problem \er{nkpp} has a unique non-negative classical solution fulfilling the following regularities
\begin{align}
    \|u(\cdot,t)\|_{L^k(\Omega)} \le C\left( \|U_0\|_{L^1(\R^n)},\|U_0\|_{L^\infty(\R^n)} \right)~~\mbox{for~any~~} 1 \le k \le \infty.
\end{align}
\end{theorem}
\noindent\textbf{Proof.} In the first place, we will prepare necessary regularities for the global existence of solutions. The derivation will be conducted step by step. Firstly, we will give the a priori estimates of $\|u(\cdot,t)\|_{L^k(\R^n)}$ for any $k>\max\left\{ \frac{(n-2)(\alpha-1)}{2},1  \right\}$. On this basis, we will show that $\|u(\cdot,t)\|_{L^k(\R^n)}$ is uniformly bounded in time for two cases of $1<k\le \alpha$ and $\alpha<k<\infty.$ Secondly, it follows the uniformly boundedness in time of solutions by making use of the iterative method. Hence we close the crucial part. Finally, combining with the standard parabolic theory we have the desired results.

\noindent{\it\textbf{Step 1.}}(A priori estimates) Multiplying \er{nkpp} by $ku^{k-1}(k \ge 1)$ we obtain
\begin{align*}
\frac{d}{dt} \int_{\R^n}  u^k dx + \frac{4(k-1)}{k} \int_{\R^n}
|\nabla u^{\frac{k}{2}} |^2 dx+ k \int_{\R^n} u dx \int_{\R^n} u^{k+\alpha-1}
dx = k M_0 \int_{\R^n} u^{k+\alpha-1} dx.
\end{align*}
Letting
\begin{align*}
w^{\frac{1}{a}}=u^{\frac{k}{2}},~~b=\frac{2(k+\alpha-1)}{k},~~a=\frac{2k'}{k},~~C_0=\frac{k-1}{k^2 M_0}
\end{align*}
in Lemma \ref{gnsyoung} for $k>\max\left( \frac{(n-2)(\alpha-1)}{2},1 \right)$ and $\max\left( \frac{n(\alpha-1)}{2},1 \right)<k'<k+\alpha-1$ we have
\begin{align*}
\int_{\R^n} u^{k+\alpha-1} dx \le \frac{k-1}{k^2 M_0}\|\nabla u^{k/2}\|_{L^2(\R^n)}^2+C(k,M_0)\|u\|_{L^{k'}(\R^n)}^\delta
\end{align*}
where
\begin{align*}
\delta=\frac{(1-\lambda)(k+\alpha-1)}{1-\frac{\lambda(k+\alpha-1)}{k}},~~\lambda=\frac{\frac{k}{2k'}-\frac{k}{2(k+\alpha-1)}}{\frac{k}{2k'}-\frac{n-2}{2n}} \in (0,1).
\end{align*}
Thus it yields
\begin{align}\label{0927}
&\frac{d}{dt}\int_{\R^n} u^k dx +k \int_{\R^n} u dx\int_{\R^n} u^{k+\alpha-1} dx+\frac{3(k-1)}{k} \|\nabla u^{k/2}\|_{L^2(\R^n)}^2 \nonumber\\
\le & C(k) \|u\|_{L^{k'}(\R^n)}^\delta.
\end{align}
Furthermore, with the help of the H\"{o}lder inequality, in case of $\max\left( \frac{n(\alpha-1)}{2},1 \right)<k'<k+\alpha-1$ one has
\begin{align}\label{09270}
\|u\|_{L^{k'}(\R^n)}^\delta & \le \|u\|_{L^{k+\alpha-1}(\R^n)}^{\delta \theta} \|u\|_{L^1(\R^n)}^{(1-\theta)\delta} \nonumber \\
& = \left( \|u\|_{L^{k+\alpha-1}(\R^n)}^{k+\alpha-1} \|u\|_{L^1(\R^n)}  \right)^{\frac{\delta \theta}{k+\alpha-1}} \|u\|_{L^1(\R^n)}^{\delta \left( 1-\theta-\frac{\theta}{k+\alpha-1} \right)}
\end{align}
where $\theta=\frac{1-\frac{1}{k'}}{1-\frac{1}{k+\alpha-1}}.$ A direct calculation shows that
\begin{align}\label{092700}
\frac{\delta \theta}{k+\alpha-1}<1
\end{align}
if and only if
\begin{align}
1 < \alpha <1+2/n.
\end{align}
Now we take $k'=\frac{k+\alpha-1+1}{2}\in (1,k+\alpha-1)$ such that
\begin{align*}
1-\theta-\frac{\theta}{k+\alpha-1}=0.
\end{align*}
Hence using the Young inequality from \er{09270} and \er{092700} we have
\begin{align*}
C(k,M_0) \|u\|_{L^{k'}(\R^n)}^\delta & \le \left( \|u\|_{L^{k+\alpha-1}(\R^n)}^{k+\alpha-1} \|u\|_{L^1(\R^n)}  \right)^{\frac{\delta \theta}{k+\alpha-1}} \\
& \le \frac{k}{2} \|u\|_{L^{k+\alpha-1}(\R^n)}^{k+\alpha-1} \|u\|_{L^1(\R^n)}+\overline{C}(k,M_0).
\end{align*}
Therefore together with \er{0927} one has that for any $k \ge \max\left\{ 1,\frac{(n-2)(\alpha-1)}{2} \right\}=1$
\begin{align}\label{092701}
\frac{d}{dt}\int_{\R^n} u^k dx +\frac{k}{2} \int_{\R^n} u dx\int_{\R^n} u^{k+\alpha-1} dx+\frac{3(k-1)}{k} \|\nabla u^{k/2}\|_{L^2(\R^n)}^2
\le \overline{C}(k,M_0).
\end{align}

\noindent{\it\textbf{Step 2.}}($L^k$ estimates for $1<k\le \alpha$)
Firstly, by the H\"{o}lder inequality and the Young inequality one has
\begin{align*}
\|u\|_{L^\alpha(\R^n)}^\alpha &\le \left( \|u\|_{L^1(\R^n)} \|u\|_{L^{2\alpha-1}(\R^n)}^{2\alpha-1} \right)^{1/2} \\
& \le \frac{\alpha}{2} \|u\|_{L^1(\R^n)} \|u\|_{L^{2\alpha-1}(\R^n)}^{2\alpha-1}+\frac{1}{2\alpha}.
\end{align*}
Hence letting $k=\alpha$ in \er{092701} one obtains
\begin{align*}
\frac{d}{dt}\int_{\R^n}u^\alpha dx + \int_{\R^n}u^\alpha dx \le C(\alpha,M_0)
\end{align*}
which assures the following uniform estimate in time
\begin{align*}
\int_{\R^n} u^\alpha dx \le \|U_0\|_{L^\alpha(\R^n)}^\alpha e^{-t}+C(\alpha,M_0) \le \|U_0\|_{L^\alpha(\R^n)}^\alpha +C(\alpha,M_0).
\end{align*}
Besides, due to
\begin{align*}
\frac{d}{dt}\int_{\R^n} u dx=\int_{\R^n}u^\alpha dx \left(M_0-\int_{\R^n} u dx \right)
\end{align*}
we have
\begin{align*}
\frac{d}{dt}\left(M_0-\int_{\R^n} u dx \right)=-\int_{\R^n}u^\alpha dx \left(M_0-\int_{\R^n} u dx \right).
\end{align*}
This leads to
\begin{align*}
\left(M_0-m_0\right) e^{-\left( \|U_0\|_{L^\alpha(\R^n)}^\alpha +C(\alpha,M_0) \right)~t } \le M_0-\int_{\R^n} u dx =\left( M_0-m_0  \right) e^{-\int_0^t \|u(s)\|_{L^\alpha(\R^n)}^\alpha ds} \le M_0-m_0
\end{align*}
where we have used the lower boundedness of $\int_{\R^n} u^\alpha dx$ and $m_0<M_0$.

\noindent{\it\textbf{Step 3.}}($L^k$ estimates for $\alpha<k<\infty$) In this step, we take advantage of $\|\nabla u^{k/2}\|_{L^2(\R^n)}^2$ to estimate $\|u\|_{L^k(\R^n)}$. Letting
\begin{align*}
w^{1/a}=u^{k/2},~~b=2,~~a=\frac{k+\alpha}{k},~~C_0=\frac{k-1}{k}
\end{align*}
in Lemma \ref{gnsyoung} follows that for $\alpha<k<\infty$
\begin{align*}
\|u\|_{L^k(\R^n)}^k & \le \frac{k-1}{k} \|\nabla u^{k/2} \|_{L^2(\R^n)}^2+C(k) \|u\|_{L^{\frac{k+\alpha}{2}}(\R^n)}^k \\
& \le \frac{k-1}{k} \|\nabla u^{k/2} \|_{L^2(\R^n)}^2+C(k) \left( \|u\|_{L^{k+\alpha-1}(\R^n)}^{k+\alpha-1} \|u\|_{L^1(\R^n)} \right)^{\frac{k}{k+\alpha}} \\
& \le \frac{k-1}{k} \|\nabla u^{k/2} \|_{L^2(\R^n)}^2+\frac{k}{2} \|u\|_{L^{k+\alpha-1}(\R^n)}^{k+\alpha-1} \|u\|_{L^1(\R^n)}+C(k)
\end{align*}
because of the H\"{o}lder inequality and the Young inequality. Substituting the above inequality into \er{092701} leads to
\begin{align*}
\frac{d}{dt} \int_{\R^n} u^k dx + \int_{\R^n} u^k dx \le C(k,M_0).
\end{align*}
Thus for any $\alpha<k<\infty$
\begin{align*}
\int_{\R^n}u^k dx \le \|U_0\|_{L^k(\R^n)}^k+C(k,M_0).
\end{align*}

\noindent{\it\textbf{Step 4.}}($L^\infty$ estimate) The route we shall follow here for the uniformly boundedness of the solution is analogous to step 4 of Theorem \ref{critical} except the exponent $\alpha$ changing from $\alpha=1+2/n$ to $\alpha<1+2/n$. Therefore, without further comment, we have
\begin{align*}
\|u\|_{L^\infty(\R^n)} \le C\left( \|U_0\|_{L^1(\R^n)},\|U_0\|_{L^\infty(\R^n)} \right).
\end{align*}
Combining with $u \in L^1 \cap L^\alpha(\R^n)$ we have that \er{nkpp} admits a global classical solution followed by the standard parabolic theory for the semilinear parabolic equation, which completes the proof. $\Box$

\section{Remarks}\label{sec5}
\def\theequation{5.\arabic{equation}}\makeatother
\setcounter{equation}{0}
\def\thetheorem{5.\arabic{theorem}}\makeatother
\setcounter{theorem}{0}

\indent(1) There is an interesting relationship between the hypotheses of Theorem \ref{critical} and Theorem \ref{subcritical}. Actually, plugging $u_\lambda(x,t)=\lambda^n u\left( \lambda x,\lambda^2 t \right)$ into \er{nkpp}, it's easy to verify that $u_\lambda(x,t)$ is also a solution of \er{nkpp} and the scaling preserves the $L^1$ norm in space, the diffusion term $\lambda^{n+2}\Delta u(\lambda x,\lambda^2 t)$ has the same scaling as the reaction term $\lambda^{n \alpha} u^\alpha \left( M_0-\int_{\R^n} u dx \right)(\lambda x,\lambda^2 t)$ if and only if
$$
\alpha=1+2/n.
$$
Firstly, observing the rescaled equation we can see that when $n \alpha<n+2$ (the subcritical case), for low density (small $\lambda$), the reaction dominates thus it prevents spreading. While for high density (large $\lambda$), the diffusion controls the reaction and thus blow-up is precluded. Therefore, in this case, the solution will exist globally (Theorem \ref{subcritical}). Secondly, for $n\alpha=n+2$ (the critical case), similar to \cite{bj07}, we guess that there is a critical value separating the global existence and the finite time blow-up, Theorem \ref{critical} provides the first step for the global existence. Thirdly, if $n \alpha>n+2$ (the supercritical case), then the diffusion manipulates for low density and the density had infinite-time spreading, the reaction dominates for high density and the solution has finite time blow-up. Hence, both the global existence and finite time blow-up occur for $n \alpha>n+2.$

(2) For $\alpha>1+2/n,$ \cite{weis81} showed that small initial data in $L^{\frac{n(\alpha-1)}{2}}(\R^n)$ can achieve global existence of the solution to \er{fujita1}. The result is the same to the following by scaling in space
\begin{align}\label{10199}
\left\{
  \begin{array}{ll}
v_t=\Delta v+M_0 v^\alpha , \quad x \in \R^n,~~t>0, \\
v|_{t=0}=v(x,0) \ge 0,\quad x \in \R^n.
  \end{array}
\right.
\end{align}
We claim that the solution of \er{10199} is a supersolution of \er{nkpp} because
\begin{align}
\frac{\partial}{\partial t} v-\Delta v=M_0 v^\alpha \ge v^\alpha \left( M_0-\int_{\R^n} v dx \right).
\end{align}
One concludes that, when $U_0(x) \le v(x,0),$ the solution of \er{nkpp} satisfies $u(x,t)\le v(x,t)$ for all times $t$ as long as the solution exists. This is sufficient to prove that the solution of \er{nkpp} exists globally for small $\|U_0(x)\|_{L^{\frac{n(\alpha-1)}{2}}(\R^n)}$. However, there is a conjecture that the dampening term $-u^\alpha \int_{\R^n} u dx$ in \er{nkpp} can relax the sufficient condition on the initial value $U_0$ for global existence. Maybe it's not necessary to impose small conditions on the initial data for global existence of solutions to \er{nkpp}.

\end{document}